\documentclass[12pt]{article}
\title{A finiteness theorem on symplectic singularities}
\author{Yoshinori Namikawa}
\date{ }
\usepackage{amsmath,amssymb,amsthm, amscd}
\chardef\bslash=`\\
\newtheorem{Thm}{Theorem}[section]
\newtheorem{Cor}[Thm]{Corollary}
\newtheorem{Lem}[Thm]{Lemma}
\newtheorem{Prop}[Thm]{Proposition}

\def\0{{\mathcal O}}

\begin{document}
\maketitle

\begin{abstract}
An affine symplectic singularity $X$ with a good $\mathbf{C}^*$-action is called a conical symplectic variety.  In this article we prove the following theorem. For fixed positive integers $N$ and $d$, there are only finite number of conical symplectic varieties of dimension $2d$ with maximal weights $N$, up to an isomorphism.   

To prove the main theorem, we first relate a conical symplectic variety with a log Fano klt pair, which has a contact structure. By the boundedness result for log Fano klt pairs with fixed Cartier index, we prove that  conical symplectic varieties of a fixed dimension and with a fixed maximal weight form a bounded family. Next we prove the rigidity of conical symplectic varieties by using Poisson deformations.
\vspace{0.2cm}

{\em 2010 Mathematics Classification}: 14D15, 14J17, 14J45, 32S30, 53D10, 53D17 


\end{abstract} 

\section{Introduction}

An affine variety $X$ is conical if $X$ can be written as $\mathrm{Spec} R$ with 
a finitely generated domain $R$ over $\mathbf{C}$ which is positively graded:  
$R = \oplus_{i \geq 0}R_i$ where $R_0 = \mathbf{C}$. The grading determines a 
$\mathbf{C}^*$-action on $X$ and the origin $0 \in X$ defined by the maximal ideal 
$m := \oplus_{i > 0}R_i$ is a unique fixed point of the $\mathbf{C}^*$-action. We 
often say that $X$ has a good $\mathbf{C}^*$-action in such a situation. 
   
A conical symplectic variety $(X, \omega)$ is a pair of a conical normal affine variety $X$ and a holomorphic symplectic 2-form $\omega$ on the regular locus $X_{reg}$, where 
(i) $\omega$ extends to a holomorphic 2-form on a resolution $Z \to X$ of $X$, and 
(ii) $\omega$ is homogeneous with respect to the $\mathbf{C}^*$-action. 

Conical symplectic varieties play an important role in algebraic geometry (cf. \cite{Be}, \cite{Na 4}) and 
geometric representation theory (cf. \cite{BPW}, \cite{BLPW}). For examples, among them are 
nilpotent orbit closures of a semisimple complex Lie algebra (cf. \cite{C-M}), Slodowy slices to 
such nilpotent orbits \cite{Slo} and Nakajima quiver varieties \cite{Nak}.  

Two conical symplectic varieties $(X_1, \omega_1)$ and $(X_2, \omega_2)$ are called 
isomorphic if there is a $\mathbf{C}^*$-equivariant isomorphism $\varphi : X_1 \to X_2$ 
such that $\omega_1 = \varphi^*\omega_2$. 

Take a set of minimal homogeneous generators $\{z_0, ..., z_n\}$ of the $\mathbf{C}$-algebra $R$. We may assume that $wt(z_0) \le wt(z_1) \le ... \le wt(z_n)$ and further assume that 
the greatest common divisor of them is $1$. We put $a_i := wt(z_i)$. Then the $n+1$-tuple  
$(a_0, ...., a_n)$ are uniquely determined by the graded algebra $R$. We call the number $a_n$ the maximal weight of $R$. Our main result is the following. 

{\bf Theorem}. {\em For positive integers $N$ and $d$, there are only finite number of 
conical symplectic varieties of dimension $2d$ with maximal weights $N$, up to an isomorphism. }\vspace{0.2cm}

{\bf Example}. We must bound the maximal weight of $X$ to have Theorem. 
In fact, $A_{2n-1}$ surface singularities 
$$ f: = x^{2n} + y^2 + z^2 = 0, \; \omega := Res(dx \wedge dy \wedge dz/f)$$ 
are conical symplectic varieties of dim $2$; their weights are $(1,n,n)$ and the 
maximal weights are not bounded above. \vspace{0.2cm}

The proof of Theorem consists of two parts. First we shall relate a conical symplectic 
variety $X$ of dimension $2d$ with a contact Fano orbifold $\mathbf{P}(X)^{orb}$ of 
dimension $2d -1$. The underlying variety $\mathbf{P}(X)$ of the orbifold is equipped with 
a divisor $\Delta$ with standard coefficients and $(\mathbf{P}(X), \Delta)$ is a  
log Fano variety with klt singularities. If we fix the maximal weight $N$ of $X$, then 
the Cartier index of $-K_{\mathbf{P}(X)} - \Delta$ is bounded above by a constant depending 
only on $d$ and $N$. By a recent result of Hacon, McKernan and Xu \cite{H-M-X}, the set 
of all such log Fano varieties forms a bounded family. This fact enables us to construct a 
flat family of conical symplectic varieties over a quasi projective base so that any conical symplectic variety of dimension $2d$ with the maximal weight $N$ appears somewhere in 
this family (Proposition \ref{(2.10)}). 

Second we shall prove that all fibres of the flat family on the same connected component 
are isomorphic as conical symplectic varieties (Proposition \ref{(3.2)}). Notice that a symplectic variety has a natural Poisson structure and the family constructed above can be regarded as a Poisson deformation of the symplectic variety. A conical symplectic variety $X$ has a universal Poisson deformation over an affine space (\cite{Na 3}). The central fibre $X$ of the universal family has a $\mathbf{C}^*$-action, but no nearby fibre does. This fact means  that $X$ is rigid under a Poisson deformation together with the $\mathbf{C}^*$-action 
(Corollary \ref{(3.2)}) and Proposition \ref{(3.3)} follows. Theorem is a corollary to Proposition \ref{(2.10)} and Proposition \ref{(3.3)}.

\section{Contact orbifolds}

In this section $(X, \omega)$ is a conical symplectic variety of dimension $2d$ with the 
maximal weight $N$. By definition $\omega$ is homogeneous with respect to the 
$\mathbf{C}^*$-action. We denote by $l$ the degree (weight) of $\omega$. 
By \cite{Na 2}, Lemma 2.2, we have $l > 0$. 
 
By using the minimal homogeneous generators in the introduction we have   
a surjection from the polynomial ring to $R$: 
$$\mathbf{C}[x_0, ..., x_n] \to R$$ which sends each $x_i$ to $z_i$.  
Correspondingly $X$ is embedded in 
$\mathbf{C}^{n+1}$.   
The quotient variety 
$\mathbf{C}^{n+1}-\{0\}/\mathbf{C}^*$ by the $\mathbf{C}^*$-action 
$(x_0, ..., x_n) \to (t^{a_0}x_0, ..., t^{a_n}x_n)$ is the weighted projective space 
$\mathbf{P}(a_0, ..., a_n)$. We put $\mathbf{P}(X) := X - \{0\}/\mathbf{C}^*$. 
By definition $\mathbf{P}(X)$ is a closed subvariety of 
$\mathbf{P}(a_0, ..., a_n)$.  Put $W_i := \{ x_i = 1\} \subset \mathbf{C}^{n+1}$. 
Then the projection map $p: \mathbf{C}^{n+1}- \{0\} \to \mathbf{P}(a_0, ..., a_n)$ 
induces a map $p_i: W_i \to \mathbf{P}(a_0, ..., a_n)$, which is a finite Galois covering 
of the image. The collection $\{p_i\}$ defines a smooth orbifold structure on 
$\mathbf{P}(a_0, ..., a_n)$ in the sense of \cite{Mum 1}, \S 2. 
More exactly, the following are satisfied 

(i) For each $i$, $W_i$ is a smooth variety, $p_i: W_i \to p_i(W_i)$ 
is a finite Galois covering\footnote{The precise definition of an orbifold only needs a 
slightly weaker condition: $p_i: W_i \to \mathbf{P}(a_0, ..., a_n)$ factorizes as 
$W_i \stackrel{q_i}\to W_i/G_i \stackrel{r_i}\to \mathbf{P}(a_0, ..., a_n)$ where $G_i$ is a 
finite group and $r_i$ is an \'{e}tale map.}, and   
$\cup \mathrm{Im}(p_i) = \mathbf{P}(a_0, ..., a_n)$.    

(ii) Let $(W_i \times_{\mathbf{P}(a_0, ..., a_n)}W_j)^n$ denote the normalization 
of the fibre product $W_i \times_{\mathbf{P}(a_0, ..., a_n)}W_j$. Then 
the maps $(W_i \times_{\mathbf{P}(a_0, ..., a_n)}W_j)^n \to W_i$ and 
$(W_i \times_{\mathbf{P}(a_0, ..., a_n)}W_j)^n \to W_j$ are both \'{e}tale maps.     
\vspace{0.2cm}

The orbifold $\mathbf{P}(a_0, ..., a_n)$ 
admits an orbifold line bundle $O_{\mathbf{P}(a_0, ..., a_n)}(1)$.  
Put $D_i := \{x_i = 0\} \subset \mathbf{P}(a_0, ..., a_n)$ and 
$D := \cup D_i$. Since $x_i$ are minimal generators, $\bar{D}:= \mathbf{P}(X) \cap D$ 
is a divisor of $\mathbf{P}(X)$.  
Let $$\bar{D}= \cup \bar{D}_{\alpha}$$ 
be the decomposition into irreducible components \footnote{The index $\alpha$ 
is usually different from the original index $i$ of $D_i$ because $D_{i_1} \cap ... 
\cap D_{i_k} \cap \mathbf{P}(X)$ may possibly become an irreducible component 
of $\bar{D}$ or $D_i \cap \mathbf{P}(X)$ may split into more than two irreducible 
components of $\bar{D}$.}.   

The map $p: X - \{0\} \to \mathbf{P}(X)$ is a $\mathbf{C}^*$-fibre bundle over 
$\mathbf{P}(X) - \bar{D}$. But a fibre over a general point of $\bar{D}_{\alpha}$ may 
possibly be a multiple fibre. We denote by $m_{\alpha}$ its multiplicity.

By putting $U_i := X \cap W_i$ and $\pi_i := p_i\vert_{U_i}$, the collection  
$\{\pi_i: U_i \to \mathbf{P}(X)\}$ of covering maps induces a (not necessarily smooth) 
orbifold structure on $\mathbf{P}(X)$.    
Namely, we have 

(i) For each $i$, $U_i$ is a normal variety and $\pi_i: U_i \to \pi_i(U_i)$ is a finite 
Galois covering. $\cup \mathrm{Im}(\pi_i) = \mathbf{P}(X)$  
(ii) The maps $(U_i \times_{\mathbf{P}(X)}U_j)^n \to U_i$ and 
$(U_i \times_{\mathbf{P}(X)}U_j)^n \to U_j$ are both \'{e}tale maps.
\vspace{0.2cm}

We put 
$\mathcal{L} := O_{\mathbf{P}(a_0, ..., a_n)}(1)\vert_{\mathbf{P}(X)}$, which is an 
orbifold line bundle on $\mathbf{P}(X)$. We call $\mathcal{L}$ the tautological 
line bundle. Then $X - \{0\} \to \mathbf{P}(X)$ can be regarded as an orbifold $\mathbf{C}^*$-bundle $(\mathcal{L}^{-1})^{\times}$.   
Notice that $X$ has only rational Gorenstein singularities and, in particular, the 
log pair $(X, 0)$ of the $X$ and the zero divisor has klt singularities. 
We define a {\bf Q}-divisor $\Delta$ by 
$$\Delta := \sum (1 - 1/m_{\alpha})\bar{D}_{\alpha}.$$ 
The following lemma (\cite{Na 1}, \S 1, Lemma)  will be a key step toward our main 
theorem.

\begin{Lem} \label{(2.1)} The pair   
$(\mathbf{P}(X), \Delta)$ is a log Fano variety, that is, $(\mathbf{P}(X), \Delta)$ 
has klt singularities and $-(K_{\mathbf{P}(X)} + \Delta)$ is an ample $\mathbf{Q}$-divisor. 
\end{Lem}

Moreover, the symplectic structure on $X$ induces a contact orbifold structure on 
$\mathbf{P}(X)$(\cite{Na 2}, Theorem 4.4.1). We shall briefly explain this. First of all, a 
contact structure on a complex manifold $Z$ of dimension $2d - 1$ is an exact 
sequence of vector bundles 
$$  0 \to E \stackrel{j}\to \Theta_Z \stackrel{\theta}\to L \to 0, $$ with a vector bundle 
$E$ of rank $2d-2$ and a line bundle $L$. Here $\theta$ induces a 
pairing map $$E \times E \to L \;\;  (x,y) \to \theta([j(x), j(y)])$$ and we require that 
it is non-degenerate. If $Z$ admits such a contact structure, then we have 
$-K_Z \cong L^{\otimes d}$.  
The map $\theta$ can be regarded as a section of $\Omega^1_Z 
\otimes L$ and we call it the contact form. Moreover, $L$ is called a contact line 
bundle. 

We can slightly generalize this notion to a singular variety $Z$. 
Let us assume that $Z$ is a normal variety of dimension $2d - 1$ and let $L$ be a
line bundle on $Z$. If $Z_{reg}$ admits a contact structure with the contact line bundle 
$L\vert_{Z_{reg}}$, then we call it a contact structure on $Z$. The twisted 1-form 
$\theta \in \Gamma (Z_{reg}, \Omega^1_Z \otimes L\vert_{Z_{reg}})$ is also called the contact form. 

We now go back to our situation. As explained above, $\mathbf{P}(X)$ admits orbifold 
charts $U_i \to \mathbf{P}(X)$. The orbifold line bundle $O_{\mathbf{P}(a_0, ..., a_n)}(1)$ 
restricts to a line bundle $L_i$ on $U_i$ and the collection $\{L_i\}$ determines an 
orbifold line bundle $\mathcal{L}$ on $\mathbf{P}(X)$. We then have a contact structure 
on each $U_i$ with the contact line bundle $L_i^{\otimes l}$, where $l = wt(\omega)$.  
Let us denote by $\theta_i$ the 
contact form. Notice that  $\theta_i$ is a  section of  $\Omega^1_{(U_i)_{reg}}\otimes L_i^{\otimes l}\vert_{(U_i)_{reg}}$. Consider the diagram 
$$U_i \stackrel{p_i}\leftarrow (U_i \times_{\mathbf{P}(X)}U_j)^n \stackrel{p_j}\to U_j.$$ 
Then we have $p_i^*(\theta_i) = p_j^*(\theta_j)$ for all $i$ and $j$.  
Thus $\theta := \{\theta_i\}$ can be regarded as a section of ${\mathcal Hom}(\Theta_{\mathbf{P}(X)^{orb}}, \mathcal{L}^{\otimes l})$. 
We call  the pair $(\theta, \mathcal{L}^{\otimes l})$ a contact orbifold structure on $\mathbf{P}(X)^{orb}$ and the orbifold line bundle $\mathcal{L}^{\otimes l}$ is called its contact line bundle. 
Similarly to the ordinary case, we have an isomorphism  
$-K_{\mathbf{P}(X)^{orb}} \cong \mathcal{L}^{\otimes ld}$ of orbifold line bundles.   
           
By the construction of $\mathbf{P}(X)$, the orbifold line bundle $\mathcal{L}^{\otimes N !}$ is 
a usual line bundle on $\mathbf{P}(X)$ and so is $-K^{\otimes N !}_{\mathbf{P}(X)^{orb}}$.  
Notice here that $K_{\mathbf{P}(X)^{orb}} = p^*(K_{\mathbf{P}(X)} + \Delta)$, where 
$p: \mathbf{P}(X)^{orb} \to \mathbf{P}(X)$ is the natural map. Therefore  
$N !(K_{\mathbf{P}(X)} + \Delta)$ is a Cartier divisor.  \vspace{0.2cm}

{\bf Theorem}(Hacon, McKernan, Xu)[H-M-X, Corollary 1.8]  {\em Let $m$ and $r$ be fixed positive integers.  Let $\mathcal{D}$ be the set of klt log Fano pairs $(Y, \Delta)$ such that $\dim Y = m$ and 
$-r(K_Y + \Delta)$ are ample Cartier divisors. Then $\mathcal{D}$ forms a bounded family.} 
\vspace{0.2cm}

In particular, the self intersection number $(-K_Y - \Delta)^m$ is bounded above by some 
constant depending only on $m$ and $r$.  

We are now going to apply the theorem above by putting  $r = N !$.

\begin{Lem} \label{(2.2)} The weight $l$ of $\omega$ is bounded above by some constant depending only on $d$ and $N$.  
\end{Lem}

{\em Proof}. Since $-K_{\mathbf{P}(X)^{orb}} \cong \mathcal{L}^{\otimes ld}$ and 
$\mathcal{L}^{\otimes N !} = p^*L$ for  $L \in \mathrm{Pic}(\mathbf{P}(X))$, we have 
$-(K_{\mathbf{P}(X)} + \Delta) \sim_{\mathbf Q} ld/N! \cdot L$. 
Here  $(-K_{\mathbf{P}(X)} - \Delta)^{2d-1}$ is bounded above by a 
constant depending on $N$ and $d$. On the other hand, $L^{2d-1}$ is a positive integer; 
this implies that $l$ must be bounded above by a constant depending on $N$ and $d$.       
Q.E.D.  

\begin{Lem} \label{(2.3)} The number of the minimal homogeneous generators of $R$ 
is bounded above by some constant depending only on $d$ and $N$. 
\end{Lem}

{\em Proof}.    
By the theorem of Hacon, McKernan and Xu, there is a positive integer $q$ (which is a 
multiple of $r$) depending 
only on $r$ and $d$ such that $q(-K_{\mathbf{P}(X)} + \Delta)$ is a very ample Cartier 
divisor and $h^0(\mathbf{P}(X), q(-K_{\mathbf{P}(X)} - \Delta) = h^0(\mathbf{P}(X), -K^{\otimes q}_{\mathbf{P}(X)^{orb}})$ is bounded above by a constant depending on $r$ and $d$. 

There are only finitely many possibilities for a weight of $R$ because the weight is smaller than 
$N$ or equal $N$. Since $r = N !$, the integer $q$ is a multiple of any possible weight. 
 
Note that $-K^{\otimes q}_{\mathbf{P}(X)^{orb}} \cong \mathcal{L}^{\otimes qld}$. 
Take an arbitrary weight, say $a$. Suppose that  
exactly $s$ elements, say, $z_1$, ..., $z_s$ have the weight $a$ among the minimal homogeneous generators. Note that these are elements of 
$H^0(\mathbf{P}(X), \mathcal{L}^{\otimes a})$. Write $q = q' a$. Then 
$(z_1)^{q'dl}$, $(z_1)^{q'dl-1}z_2$, ...., $(z_1)^{q'dl-1}z_s$ are linearly independent elements of 
$H^0(\mathbf{P}(X), \mathcal{L}^{\otimes qld})$. In fact, suppose to the contrary that
there is a non-trivial relation $\lambda_1(z_1)^{q'dl} + \lambda_2(z_1)^{q'dl-1}z_2 + ... +  \lambda_s(z_1)^{q'dl-1}z_s = 0.$  Then we have an equality 
$$(z_1)^{q'dl-1}\cdot (\lambda_1z_1 + ... + \lambda_s z_s) = 0$$ in $R = \oplus_{i \geq 0} H^0(\mathbf{P}(X), \mathcal{L}^{\otimes i})$. Since $R$ is a domain, we conclude that 
$z_1 = 0$ or $\Sigma \lambda_i z_i = 0$. But, by the assumption, both $z_1$ and 
$\Sigma \lambda_i z_i$ are nonzero, which is absurd; hence 
$(z_1)^{q'dl}$, $(z_1)^{q'dl-1}z_2$, ...., $(z_1)^{q'dl-1}z_s$ are linearly independent. 
This means that 
$$s \le h^0(\mathbf{P}(X), \mathcal{L}^{\otimes qld}) =  h^0(\mathbf{P}(X), -K^{\otimes q}_{\mathbf{P}(X)^{orb}}).$$ In particular, $s$ is bounded above by a constant depending only on 
$d$ and $N$.  Q.E.D.       
\vspace{0.2cm}

The graded coordinate ring of a weighted projective space is called a {\em weighted polynomial ring}.

\begin{Cor} \label{(2.4)}  Fix positive integers $d$ and $N$. Then there are finitely many 
weighted polynomial rings $S_1$, ..., $S_k$ such that any graded coordinate ring $R$ of a  conical symplectic variety of dimension $2d$ with the maximal weight $N$ can be realized as a quotient of some $S_i$. 
\end{Cor}

In the corollary, we put $\mathbf{P}_i := \mathrm{Proj}(S_i)$. Then we have:  

\begin{Cor} \label{(2.5)} There are flat families of closed subschemes of $\mathbf{P}_i$( $1 \le  i \le k)$: $\mathcal{Y}_i \subset \mathbf{P}_i \times T_i$ parametrized by reduced quasi projective schemes $T_i$ such that, for any conical symplectic variety $X$ of dimension 
$2d$ with the maximal weight $N$, there is a point $t \in T_i$ for some $i$ and $\mathbf{P}(X) = \mathcal{Y}_{i,t}$. 
\end{Cor}
\vspace{0.2cm}

{\em Proof}. Let $q$ be the least common multiple of all weights of the minimal homogeneous generators of $S_i$'s. Then $O_{\mathbf{P}_i}(q)$ is an ample line bundle for every $i$. Take a conical 
symplectic symplectic variety $X$ of dimension $2d$ with the maximal weight $N$. 
Then $\mathbf{P}(X)$ can be embedded in some $\mathbf{P}_i$. By the theorem of Hacon, 
McKernan and Xu, there are only finitely many possibilities of the Hilbert polynomial 
$\chi (\mathbf{P}(X), O_{\mathbf{P}(X)}(qn))$. Such closed subschemes of $\mathbf{P}_i$ 
form a bounded family. Q.E.D. 
\vspace{0.2cm}

Let $\mathcal{Y} \subset \mathbf{P} \times T$ be one of the flat families in \ref{(2.5)}.  Define a map  $f: \mathcal{Y} \to T$ to be the composite $\mathcal{Y} \to 
\mathbf{P} \times T \stackrel{pr_2}\to T$.  
Let $\{W_i \to \mathbf{P}\}$ be the orbifold charts for the weighted projective space
$\mathbf{P}$ constructed in the beginning of this section. Denote by $G_i$ the Galois group 
for $W_i \to \mathbf{P}$. Then the collection $\{W_i \times T \to \mathbf{P} \times T\}$ also gives relative orbifold charts for $\mathbf{P} \times T/T$.    
By pulling back these charts by the inclusion map $\mathcal{Y} \to \mathbf{P} \times T$, we 
have relative orbifold charts $\{\mathcal{U}_i \stackrel{\pi_i}\to \mathcal{Y}\}$ for  
$\mathcal{Y}/T$. 
If necessary, we stratify 
$T$ into a disjoint union of finite number of locally closed sets by using the generic flatness property (cf. \cite{Mum 3}, Lecture 8, Proposition) repeatedly so that all $\mathcal{U}_i$ are flat over each stratum and replace $T$ by the disjoint union of such subsets.  
Thus we may assume that $\mathcal{U}_i$ are all flat over $T$.
Let $O^{orb}_{\mathbf P}(1) := \{O_{W_i}(1)\}$ be the tautological orbifold 
line bundle on $\mathbf{P}$. Denote simply by $O_{W_i \times T}(1)$ the pullback of 
$O_{W_i}(1)$ by the projection $W_i \times T \to W_i$. 
Then $\{O_{W_i \times T}(1)\vert_{{\mathcal U}_i}\}$ gives a relative tautological orbifold line 
bundle $O^{orb}_{\mathcal Y}(1)$ on $\mathcal{Y}$. For each $j \in \mathbf{Z}$, we define 
a usual sheaf $\mathcal{L}^j$ on $\mathcal{Y}$ by  
$\mathcal{L}^j := \{\pi^{G_i}_*O_{W_i \times T}(j)\vert_{{\mathcal U}_i}\}$.    
Notice that $\mathcal{L}^j$ is flat over $T$.  
On the other hand, for $t \in T$, one can consider the orbifold structure on 
$\mathcal{Y}_t$ induced by the embedding $\mathcal{Y}_t \subset \mathbf{P}$. We similarly 
define a tautological orbifold line bundle $O^{orb}_{\mathcal{Y}_t}(1)$ and the usual sheaves 
$\mathcal{L}_t^j$ on $\mathcal{Y}_t$. We have $\mathcal{L}^j \otimes_{O_{\mathcal Y}}
O_{\mathcal{Y}_t} \cong  \mathcal{L}_t^j$. 

We define 
$$ T' := \{t \in T; \; \; f_*\mathcal{L}^j \otimes_{O_T} k(t) \cong 
H^0(\mathcal{Y}_t, \mathcal{L}_t^j) \;\; \mathrm{and} \; f_*\mathcal{L}^j \; \mathrm{are} \; \mathrm{locally} \; \mathrm{free} \; \mathrm{at} \; t \;  
 \mathrm{for} \; \mathrm{all} \; j \geq 0\}.$$   
 
\begin{Lem} \label{(2.6)} The set $T'$ is a non-empty Zariski open subset of $T$.  
\end{Lem}

{\em Proof}. First we show that there is a positive integer $j_0$ such that 
$H^1(\mathcal{Y}_t, \mathcal{L}^j_t) = 0$ for all $j \geq j_0$ and for all $t \in T$. 
Take a positive integer $q$ so that $O_{\mathbf{P}}(q)$ is a very ample line bundle on 
$\mathbf{P}$. Notice that $\mathcal{L}^j \otimes O_{\mathbf{P}}(q) \cong 
\mathcal{L}^{j+q}$ for all $j$. 
We consider the sheaves $\mathcal{L}^j$ for $j$ with $0 \le j < q$. Notice that 
they are flat over $T$.
%Since each $\mathcal{L}^j$ is flat over $T$, the Hilbert polynomials $\chi (\mathcal{Y}_t, 
%\mathcal{L}^j_t(qn))$ depend only on connected components of $T$. Since 
%$T$ is quasi projective, $T$ has only a finite number of connected components.
 
We shall prove that there is a positive integer $n_j$ such that $H^1(\mathcal{Y}_t, \mathcal{L}^j_t(qn)) = 0$ for all $n \geq n_j$ and for all $t$.
As $T$ is of finite type over $\mathbf{C}$, we can take a positive integer $n_j$ so that 
$R^pf_*\mathcal{L}^j(qn) = 0$ for all $p > 0$ and for all $n \geq n_j$ (Serre vanishing theorem). 
Notice that, for any point $t \in T$, one has $H^p(\mathcal{Y}_t, \mathcal{L}_t^j(qn)) = 0$ for $p > 2d - 1$. Fix an integer $n$ with $n \geq n_j$ and let us consider the base change map $\varphi^p(t): R^pf_*\mathcal{L}^j(qn) \otimes k(t) \to H^p(\mathcal{Y}_t, \mathcal{L}_t^j(qn))$.
By the base change theorem, when $\varphi^p(t)$ is surjective, the map $\varphi^{p-1}(t)$ is a 
surjection if and only if $R^pf_*\mathcal{L}^j(qn)$ is locally free at $t$. 
 First of all, $\varphi^{2d}(t)$ is a surjection because $H^{2d}(\mathcal{Y}_t, \mathcal{L}_t^j(qn)) 
= 0$. Since $R^{2d}f_*\mathcal{L}^j(qn) = 0$, the map $\varphi^{2d-1}(t)$ is also a surjection. 
We can repeat this argument to conclude that $\varphi^1(t)$ is a surjection because 
$R^pf_*\mathcal{L}^j(qn) = 0$ for all $p > 0$. Since $R^1f_*\mathcal{L}^j(qn) = 0$, we finally 
obtain that $H^1(\mathcal{Y}_t, \mathcal{L}^j_t(qn)) = 0$. 

%We embed $\mathbf{P}$ into a projective space $\mathbf{P}^N$ by 
%$O_{\mathbf P}(q)$. Fix $j$ so that $0 \le j < q$. By Serre's theorem there is a positive integer $m$ such that we have a surjection  $O_{\mathbf{P}^N}^{\oplus d} \to O_{\mathbf P}(j + mq)$ 
%for some $d$. For $t \in T$, there is a surjection $O_{\mathbf P}(j + mq) 
%\to \mathcal{L}^j_t(mq)$. Let $\mathcal{F}_t$ be the kernel of the composition map of 
%$O_{\mathbf{P}^N}^{\oplus d} \to O_{\mathbf P}(j + mq)$ and $O_{\mathbf P}(j + mq) 
%\to \mathcal{L}^j_t(mq)$. Then we have an exact sequence 
%$$ 0 \to \mathcal{F}_t \to O_{\mathbf{P}^N}^{\oplus d} \to \mathcal{L}^j_t(mq) \to 0.$$ 
%The Hilbert polynomial $\chi (\mathbf{P}^N, \mathcal{F}_t(n))$ only depends on the 
%connected component of $T$ that contains $t$. As $T$ has only finitely many connected 
%components, by [Mum 3, Lecture 14] and [FAG, Theorem 5.3], we have a positive constant $m_j$ such that 
%$H^2(\mathcal{F}_t(n)) = 0$ for all $n \geq m_j$ and for all 
%$t \in T$. By the exact sequence we have $H^1(\mathcal{L}^j_t((m + n)q)) = 0$ for 
%all $n \geq m_j$ and for all $t \in T$. If we put $n_j := m + m_j$, then  
%$H^1(\mathcal{Y}_t, \mathcal{L}^j_t(qn)) = 0$ for all $n \geq n_j$ and for all $t$. 

Put $\nu := \mathrm{max}\{n_0, ..., n_{q-1}\}$. Then we have 
$H^1({\mathcal Y}_t,  \mathcal{L}^j_t) = 0$ for all $j \geq q\nu$ and for all $t \in T$. 
By the base change theorem, $f_*\mathcal{L}^j$ are locally free and 
$f_*\mathcal{L}^j \otimes_{O_T} k(t) \cong H^0(\mathcal{Y}_t, \mathcal{L}_t^j)$ for 
all $j \geq q\nu$ and for all $t \in T$. 

We next consider the sheaves $\mathcal{L}^j$ for $j < q\nu$. For each such $j$, 
it is an open condition for $T$ that $f_*\mathcal{L}^j$ is locally free at $t$ and  $f_*\mathcal{L}^j \otimes_{O_T} k(t) \cong H^0(\mathcal{Y}_t, \mathcal{L}_t^j)$ 
holds. Therefore $T'$ is a non-empty Zariski open subset of $T$.  Q.E.D.  
\vspace{0.2cm}

Define 
$$ \mathcal{X} := (\mathbf{Spec}_{T} \oplus_{j \geq 0} f_*\mathcal{L}^j)  \times_T T' .$$ 

As each direct summand  $f_*\mathcal{L}^j$ is flat over $T$, the map $\mathcal{X} \to T'$ 
is flat. Let us return to \ref{(2.5)}. The construction above enables us to make a 
flat family $\mathcal{X}_i \to T'_i$ of affine schemes with good $\mathbf{C}^*$-actions 
on an open subset $T'_i$ of each $T_i$.

\begin{Cor} \label{(2.7)} There is a flat family of affine schemes with good $\mathbf{C}^*$-actions: $\mathcal{X} \to T$ parametrized by reduced quasi projective schemes $T$ such that, for any conical symplectic variety $X$ of dimension 
$2d$ with the maximal weight $N$, there is a point $t \in T$ and $X \cong  \mathcal{X}_t$ as a $\mathbf{C}^*$-variety. 
\end{Cor}

{\em Proof}. The family $\mathcal{X} \to T$ is nothing but the disjoint union of 
$\{\mathcal{X}_i \to T'_i\}$. Let $X$ be a conical symplectic variety of dim $2d$ with the 
maximal weight $N$. By \ref{(2.5)} there is a point $t$ of some $T_i$ and $\mathbf{P}(X) = \mathcal{Y}_{i,t} \subset \mathbf{P}_i$. Since the coordinate ring $R$ of $X$ is normal, 
the natural maps $H^0(\mathbf{P}_i, O_{{\mathbf P}_i}(j)) \to 
H^0(\mathbf{P}(X), O_{{\mathbf P}(X)}(j))$ are surjective for all $j \geq 0$. 
This fact implies that $t \in T'_i$ and $\mathcal{X}_{i,t} = X$. Q.E.D. 
\vspace{0.2cm}

Let $f: \mathcal{X} \to T$ be the flat family in \ref{(2.7)}. By Elkik \cite{E}, Theoreme 4, 
the set $$\mathcal{X}^0 := \{x \in \mathcal{X}; \; \mathcal{X}_{f(x)} \; \mathrm{has} \; 
\mathrm{rational} \;  \mathrm{singularities} \; \mathrm{at} \; x \}$$ is a Zariski open 
subset of $\mathcal{X}$. By the $\mathbf{C}^*$-action we also see that 
$$T^0 := \{t \in T; \mathcal{X}_t \; \mathrm{has} \; \mathrm{rational} \; 
\mathrm{singularities} \}$$ is a Zariski open subset of $T$. We take a resolution $T_{rat}$ 
of $T^0$. Note that $T_{rat}$ is the disjoint union of finitely many nonsingular quasi projective 
varieties.  We put $\mathcal{X}_{rat} := 
\mathcal{X} \times_T T_{rat}$ and let $f_{rat} :\mathcal{X}_{rat} \to T_{rat}$ be the induced 
flat family. Again by \cite{E}, $\mathcal{X}_{rat}$ has only rational singularities. In particular, it 
is normal. Notice that any conical symplectic variety of dim $2d$ with the maximal weight $N$ is realized as a fibre of this family.    

We next stratify $T_{rat}$ into the disjoint union of locally closed smooth subsets $T_{rat,i}$ 
so that $\mathcal{X}_{rat} \times_{T_{rat}}T_{rat, i} \to T_{rat,i}$ have $\mathbf{C}^*$-equivariant simultaneos resolutions. 
To obtain such a stratification, we first take a $\mathbf{C}^*$-equivariant resolution $\tilde{\mathcal{X}}_{rat} \to 
\mathcal{X}_{rat}$. By Bertini's theorem there is an open subset $T^0_{rat}$ of $T_{rat}$ such that this resolution gives  simultaneous resolutions of fibres over $T^0_{rat}$. Next stratify 
the complement $T_{rat} - T^0_{rat}$ into locally closed smooth subsets, take maximal strata 
and repeat the same for the families over them. Thus we have proved the following.

\begin{Prop} \label{(2.8)} There is a flat family of affine varieties with good $\mathbf{C}^*$-actions: $\mathcal{X} \to T$ parametrized by the disjoint union 
$T$ of a finite number of quasi projective nonsingular varieties such that  

(i) $\mathcal{X}_t$ have only rational singularities for all $t \in T$,

(ii) there is a $\mathbf{C}^*$-equivariant simultaneous resolution $\mathcal{Z} \to \mathcal{X}$ of $\mathcal{X}/T$; namely, $\mathcal{Z}_t \to \mathcal{X}_t$ are resolutions for all $t \in T$, and

(iii) for any conical symplectic variety $X$ of dimension $2d$ with the maximal weight $N$, there is a point $t \in T$ and $X \cong  \mathcal{X}_t$ as a $\mathbf{C}^*$-variety. 
\end{Prop}

Let $\mathcal{X} \stackrel{f}\to T$ and $\mathcal{Z} \stackrel{g}\to T$ be the families in \ref{(2.8)}. Let us consider the relative dualizing sheaf $\omega_{\mathcal{X}/T}$ 
of $f$. For $t \in T$, we have $\omega_{\mathcal{X}/T} \otimes_{O_{\mathcal X}} O_{\mathcal{X}_t} \cong \omega_{\mathcal{X}_t}$. The locus $T_{gor} \subset T$ where $\omega_{\mathcal{X}_t}$ is invertible is an open subset of $T$. We put $\mathcal{X}_{gor} 
:= \mathcal{X} \times_T T_{gor}$ and $\mathcal{Z}_{gor} := \mathcal{Z} \times_T T_{gor}$. 
Then $f$ and $g$ induce respectively maps $f_{gor}: \mathcal{X}_{gor} \to T_{gor}$ and 
$g_{gor}: \mathcal{Z}_{gor} \to T_{gor}$. Notice that any conical symplectic variety $X$ of 
dim $2d$ with the maximal weight $N$ still appears in some fibre of $f_{gor}$.

\begin{Prop} \label{(2.9)} (Base change theorem): Let $h: W \to S$ be a morphism of 
quasi projective schemes over $\mathbf{C}$. Assume that $W$ is normal and $\mathbf{C}^*$ 
acts on $W$ fibrewisely with respect to $h$. Let $F$ be a $\mathbf{C}^*$-linearized coherent 
$O_W$-module on $W$, which is flat over $S$. 
Then the higher direct image sheaves $R^ih_*F$ are naturally graded: $R^ih_*F = \oplus_{j \in \mathbf{Z}} (R^ih_*F)(j)$. Assume that $(R^ih_*F)(j_0)$ $(i \geq 0)$ are all coherent sheaves on $S$ for $j_0$. Then the following hold. 
 
%Then, for $s \in S$, 
%the natural homomorphism $\phi^i_s : (R^ih_*F) (j_0) \otimes_{O_T}k(s) \to 
%H^i(W_s, F_s)(j_0)$ have the following properties.} 

(a) For each $i$, the function $S \to \mathbf{Z}$ defined by 
$$s \to \dim H^i(W_s, F_s)(j_0)$$ is upper-semicontinuous on $S$.

%{\em if $\phi^i_s$ is surjective, then there is an open neighborhood $U$ of $s \in S$ such that %$\phi^i_{s'}$ are surjective for all $s' \in U$, and} 

(b) Assume that $S$ is reduced and connected. If the function $s \to \dim H^i(W_s, F_s)(j_0)$ is constant, then $(R^ih_*F)(j_0)$ is locally free sheaf on $S$ and, for all $s \in S$, 
the natural map $\phi^i_s : (R^ih_*F) (j_0) \otimes_{O_T}k(s) \to 
H^i(W_s, F_s)(j_0)$ is an isomorphism.

%{\em with the same assumption in (a), $(R^ih_*F)(j_0)$ is a free $O_S$-module in an open 
%neighborhood of $s$ if and only if $\phi^{i-1}_s$ is also surjective.}       
\end{Prop}

We can take $\mathbf{C}^*$-equivariant affine open coverings of $W$ by the theorem of 
Sumihiro (cf. \cite{K-K-M-S}, Chapter I, \S 2.).  Then the proof of \ref{(2.9)} is similar to \cite{Mum 2}, II, 5. 
\vspace{0.2cm}  

We apply this proposition to $g_{gor}: \mathcal{Z}_{gor} \to T_{gor}$ and $\Omega^k_{\mathcal{Z}_{gor}/T_{gor}}$. 
Notice that $(R^i(g_{gor})_*\Omega^k_{\mathcal{Z}_{gor}/T_{gor}})(l)$ are all coherent sheaves on $T_{gor}$ for any $l$. Let us consider the relative differential map $$((g_{gor})_*\Omega^2_{\mathcal{Z}_{gor}/T_{gor}})(l) \stackrel{d}\to ((g_{gor})_*\Omega^3_{\mathcal{Z}_{gor}/T_{gor}})(l)$$ and put $$\mathcal{F} := \mathrm{Ker}(d), \; \mathcal{G} := \mathrm{Coker}(d) $$
 
Fix an integer $l$. Then one can find a non-empty Zariski open dense subset $T_l$ of $T_{gor}$ so that, if 
$t \in T_l$, then  
both $\mathcal{F}$ and $\mathcal{G}$ are free at $t$,  $$((g_{gor})_*\Omega^2_{\mathcal{Z}_{gor}/T_{gor}})(l) \otimes k(t) \cong 
H^0(\mathcal{Z}_{gor,t}, \Omega^2_{\mathcal {Z}_{gor, t}})(l),$$ and  
$$((g_{gor})_*\Omega^3_{\mathcal{Z}_{gor}/T_{gor}})(l) \otimes k(t) \cong 
H^0(\mathcal{Z}_{gor,t}, \Omega^3_{\mathcal {Z}_{gor, t}})(l).$$ 
%and 
%$$((g_{gor})*\Omega^{2d}_{\mathcal{Z}_{gor}/T})(dl) \otimes k(t) \cong 
%H^0(\mathcal{Z}_{gor,t}, \Omega^{2d}_{\mathcal {Z}_{gor, t}})(l).$$   

A d-closed $2$-form $\omega_0$ of weight $l$ on a fibre $\mathcal{Z}_{gor,t}$ $(t \in T_l)$ is an element of $\mathrm{Ker}[H^0(\mathcal{Z}_{gor,t}, \Omega^2_{\mathcal {Z}_{gor, t}})(l)  
\stackrel{d}\to H^0(\mathcal{Z}_{gor,t}, \Omega^3_{\mathcal {Z}_{gor, t}})(l)]$. By the exact sequence 
$$ \mathcal{F}\otimes k(t) \to   
H^0(\mathcal{Z}_{gor,t}, \Omega^2_{\mathcal {Z}_{gor, t}})(l) 
\stackrel{d}\to H^0(\mathcal{Z}_{gor,t}, \Omega^3_{\mathcal {Z}_{gor, t}})(l),$$ 
$\omega_0$ comes from an element $\omega'_0 \in \mathcal{F}\otimes k(t)$. 
Then $\omega'_0$ lifts to a local section $\omega$ of $\mathcal{F}$. If we regard $\omega$ 
as a local section of $((g_{gor})_*\Omega^2_{\mathcal{Z}_{gor}/T_{gor}})(l)$, it is a d-closed 
relative 2-form extending the original $\omega_0$. 
   
Assume that $\mathcal{X}_t$ $(t \in T_l)$ is a conical symplectic variety and 
$\omega_0$ is the extension of the symplectic 2-form on $\mathcal{X}_{t, reg}$ to 
the resolution $\mathcal{Z}_t$. The wedge product $\wedge^{d}\omega_0$ is regarded as  
a section of the dualizing sheaf $\omega_{\mathcal{X}_{gor,t}}$ by the identification 
$H^0(\mathcal{Z}_{gor,t}, \Omega^{2d}_{\mathcal{Z}_{gor,t}}) \cong H^0(\mathcal{X}_{gor,t}, \omega_{\mathcal{X}_{gor,t}})$. 
Then $\wedge^{d}\omega_0$ generates the invertible sheaf $\omega_{\mathcal{X}_{gor,t}}$. 
We also see that $\wedge^{d}\omega$ generates $\omega_{\mathcal{X}_{gor}/T_{gor}}$ on 
near fibres of $\mathcal{X}_{gor,t}$.  

The argument here shows that 
$$ T_{symp,l} := \{t \in T_l; \; \mathcal{X}_{gor,t}\; \mathrm{is}\;  \mathrm{a}\; \mathrm  {conical}\;  \mathrm{symplectic}\;  \mathrm{variety}\;  \mathrm{with}\;  \mathrm{a}\;  
\mathrm{symplectic}\; \mathrm{form}\; \mathrm{of}\; \mathrm{weight}\; l\} $$ 
is an open subset of $T_l$. We have fixed an integer $l$. But, notice that the choice of 
such an $l$ is finite by \ref{(2.2)}.  

We put $\mathcal{X}_{symp, l} := \mathcal{X} \times_T T_{symp,l}$ and 
$\mathcal{Z}_{symp,l} := \mathcal{Z} \times_T T_{symp,l}$. Then  
$\mathcal{X}_{symp, l} \to T_{symp,l}$ is 
a flat family of conical symplectic varieties with symplectic forms of weight $l$ and  
$\mathcal{Z}_{symp,l} \to T_{symp,l}$ is its simultaneous resolution. 

To specify the symplectic form on each fibre of $\mathcal{X}_{symp, l} \to T_{symp,l}$, we 
consider the vector bundle $p: V_l := \mathbf{V}(\mathcal{F}^*\vert_{T_{symp,l}}) \to T_{symp,l}$. 
Each point $v \in V_l$ corresponds to a d-closed holomorphic 2-form $\omega_v$ on the regular part of $(\mathcal{X}_{symp,l})_{p(v)}$. Moreover $\omega_v$ extends to a holomorphic 2-form on the resolution $(\mathcal{Z}_{symp,l})_{p(v)}$. 
Let $V_l^0$ be the non-empty Zariski open subset of $V_l$ where $\omega_v$ is nondegenerate.  
Take the base change $\mathcal{X}_{symp,l} \times_{T_{symp.l}} V_l^0 \to V_l^0$.  
Then the regular locus of a fibre of this family is naturally equipped with a symplectic 2-form of weight $l$.  
 
Stratify $T \setminus T_l$ into locally closed smooth subsets, take maximal 
strata and repeat the same for the families over them. Then we get: 

\begin{Prop} \label{(2.10)} There is a flat family of the pairs of affine symplectic varieties with good $\mathbf{C}^*$-actions and symplectic forms: $(\mathcal{X}, \omega_{\mathcal{X}/T}) \to T$ parametrized by the disjoint union $T$ of a finite number of quasi projective nonsingular varieties such that

(i) for each connected component $T_i$ of $T$, all fibres $(\mathcal{X}_t, \omega_t)$ over 
$t \in T_i$ are conical symplectic varieties admitting symplectic forms of a fixed 
weight $l_i > 0$

(ii) there is a $\mathbf{C}^*$-equivariant simultaneous resolution $\mathcal{Z} \to \mathcal{X}$ of $\mathcal{X}/T$; namely, $\mathcal{Z}_t \to \mathcal{X}_t$ are resolutions for all $t \in T$, and 

(iii) for any conical symplectic variety $(X, \omega)$ of dimension $2d$ with the maximal weight $N$, there is a point $t \in T$ and $(X, \omega) \cong  (\mathcal{X}_t, \omega_t)$ as a $\mathbf{C}^*$-symplectic variety. 
\end{Prop}

\section{Rigidity of conical symplectic varieties}

Let $(X, \omega)$ be a conical symplectic variety with a symplectic form $\omega$ of 
weight $l$. The symplectic form $\omega$ determines a Poisson structure on $X_{reg}$. 
By the normality of $X$, this Poisson structure uniquely extends to a Poisson structure 
$\{\; , \;\}$ on $X$. Here a Poisson structure on $X$ exactly means a skew-symmetric $\mathbf{C}$-bilinear map $\{\; , \;\} : O_X \times O_X \to O_X$ which is a biderivation with respect to the 1-st and the 2-nd factors, and satisfies the Jabobi identity. We will consider a Poisson deformation of the Poisson variety. A $T$-scheme $\mathcal{X} \to T$ is called a Poisson $T$-scheme if  there is a $O_T$-bilinear Poisson bracket $\{\; , \;\}_{\mathcal{X}} : O_{\mathcal X} \times O_{\mathcal X} \to O_{\mathcal X}$, which is a biderivation, and satisfies the Jacobi identity. Let $T$ be a scheme over $\mathbf{C}$ and let $0 \in T$ be a closed point.   

A Poisson deformation of the Poisson variety $X$ over $T$ is a Poisson $T$-scheme  
$f: \mathcal{X} \to T$ together with an isomorphism $\varphi: \mathcal{X}_0 \cong X$ 
which satisfies the following conditions 

(i) $f$ is a flat surjective morphism, and 

(ii)$\{\; , \;\}_{\mathcal X}$ restricts to the original Poisson structure $\{\; , \;\}$ on 
$X$ via the identification $\varphi$.

Two Poisson deformations $(\mathcal{X}/T, \varphi)$ and $(\mathcal{X}'/T, \varphi')$ 
with the same base are equivalent if there is a $T$-isomorphism $\mathcal{X} \cong 
\mathcal{X}'$ of Poisson schemes such that it induces the identity on the central fibre. 
For a local Artinian $\mathbf{C}$-algebra $A$ with residue field $\mathbf{C}$, we define 
$\mathrm{PD}_X(A)$ to be the set of equivalence classes of Poisson deformations of $X$ over 
$\mathrm{Spec}(A)$. Then it defines a functor 
$$\mathrm{PD}_X: (Art)_{\mathbf C} \to (Set)$$ from the category of local Artinian 
$\mathbf{C}$-algebra with residue field $\mathbf{C}$ to the category of sets. 
\vspace{0.2cm}

\begin{Thm}(\cite{Na 3}, Theorem 5.5) There is a Poisson deformation $\mathcal{X}_{univ} \to \mathbf{A}^m$ of $X$ over an affine space $\mathbf{A}^m$ with $\mathcal{X}_{univ,0} = X$. This Poisson deformation has the following properties and it 
is called the universal Poisson deformation of $X$.

(i) For any Poisson deformation $\mathcal{X} \to T$ of $X$ over $T = \mathrm{Spec}(A)$ with $A \in (Art)_{\mathbf C}$, there is a unique morphism $\phi: T \to \mathbf{A}^m$ which sends the closed point of $T$ to the center $0 \in \mathbf{A}^m$ such that 
$\mathcal{X}/T$ and $\mathcal{X}_{univ} \times_{\mathbf{A}^m} T/ T$ are equivalent as Poisson deformations of $X$. 

(ii) There are natural $\mathbf{C}^*$-actions on $\mathcal{X}_{univ}$ and $\mathbf{A}^m$ induced from the $\mathbf{C}^*$-action on $X$ such that the map $\mathcal{X}_{univ} \to \mathbf{A}^m$ is $\mathbf{C}^*$-equivariant. Moreover the coordinate ring $\mathbf{C}[y_1, ..., y_m]$ of $\mathbf{A}^m$ is positively graded so that $wt(y_i) > 0$ for all $i$. 
\end{Thm}

\begin{Cor} \label{(3.2)} Let $(X, \omega)$ be a conical symplectic variety and 
let $T := \mathrm{Spec}(A)$ be a nonsingular affine curve with a base point $0 \in T$.  
Assume that $\mathcal{X} \to T$ is a Poisson deformation of $X$. Assume that  
$\mathbf{C}^*$ acts on $\mathcal{X}$ in such a way that 

(i) it induces a $\mathbf{C}^*$-action on each fibre of $\mathcal{X}/T$ and the 
$\mathbf{C}^*$-action on the central fibre coincides with the original $\mathbf{C}^*$-action 
on $X$, and 

(ii) the Poisson bracket on each fibre is homogeneous 
with respect to this action. 
 
Then there is a $\mathbf{C}^*$-equivariant Poisson isomorphism 
$f: \mathcal{X} \times_T \hat{T} \cong  X \times \hat{T}$ over $\hat{T} := \mathrm{Spec}(\hat{A})$, 
where $\hat{A}$ is the completion of $A$ along the defining ideal $m$ of $0$.     
\end{Cor}

{\em Proof}. Notice that $\hat{A} \cong \mathbf{C}[[t]]$. Put $T_n := \mathrm{Spec}
(\mathbf{C}[[t]]/(t^{n+1})$ and $X_n := \mathcal{X} \times_T T_n$. The formal 
Poisson deformation $\{X_n \to T_n\}$ determines a morphism 
$\phi: \mathrm{Spec}(\mathbf{C}[[t]]) \to \mathbf{A}^m$ by the previous theorem. 
By the assumption (i), $\mathrm{Im}(\phi)$ is contained in the $\mathbf{C}^*$-fixed locus 
of $\mathbf{A}^m$. By the last property in (ii) of Theorem, this means that 
$\phi$ is the constant map to the origin of $\mathbf{A}^m$.   
We are now going to construct an isomorphisms between formal Poisson deformations $\{X_n\}_{n \geq 0} \stackrel{\beta_n}\to \{X \times T_n\}_{n \geq 0}$, where the right hand side is a trivial Poisson deformation of $X$. Assume that we already have $\beta_{n-1}$. 
Since $\phi$ is constant, we have an equivalence $X_n \stackrel{\beta'_n}\cong X \times 
T_n$ of Poisson deformations of $X$. We put $\beta'_{n-1} := \beta'_n \vert_{X_{n-1}}$. 
Then $\gamma_{n-1} := \beta^{-1}_{n-1} \circ \beta'_{n-1}$ is a Poisson automorphism of $X_{n-1}$. By (the proof of) Corollary 2.5 of \cite{Na 3}, $\gamma_{n-1}$ lifts to a Poisson 
automorphism $\gamma_n$ of $X_n$. Define $\beta_n := \beta'_n \circ \gamma_n^{-1}$. 
Then $\beta_n$ is a Poisson isomorphism from $X_n$ to $X \times T_n$ extending 
$\beta_{n-1}$.     

Note that $\{X \times T_n\}_{n \geq 0}$ has a natural $\mathbf{C}^*$-action induced by 
the $\mathbf{C}^*$-action of $X$. 
By the isomorphisms $\{\beta_n\}_{n \geq 0}$ above, this $\mathbf{C}^*$-action induces a 
$\mathbf{C}^*$-action on $\{X_n\}_{n \geq 0}$. 
On the other hand, $\{X_n\}_{n \geq 0}$ has a $\mathbf{C}^*$-action inherited from 
$\mathcal{X}$. We will construct a Poisson automorphism $\{\psi_n\}_{n \geq 0}$ 
of $\{X_n\}_{n \geq 0}$ inductively so that these two $\mathbf{C}^*$-actions are 
compatible. At first we put $\psi_0 := id$. Assume that we are given a Poisson automorphism 
$\psi_n$ which makes two $\mathbf{C}^*$-actions compatible. 
By (the proof of) Corollary 2.5 of \cite{Na 3} $\psi_n$ lifts to a Poisson automorphism $\psi'_{n+1}$ of $X_{n+1}$. We denote by $\zeta_1$, $\zeta_2 \in \Gamma (X, \Theta_{X_{n+1}/T_{n+1}})$ respectively the relative vector fields generating the 1-st and 2-nd $\mathbf{C}^*$-actions. 
Let $\zeta \in \Gamma (X, \Theta_X)$ be the vector field (Euler vector field) generating the $\mathbf{C}^*$-action. Notice that $\zeta_1\vert_X = \zeta_2\vert_X = \zeta$. 
We write $$(\psi'_{n+1})_*\zeta_1 - \zeta_2 = t^{n+1}\cdot \Sigma v_i, $$ with 
$v_i \in \Gamma (X, \Theta_{X})(i)$. In other words, $v_i$ is a homogeneous vector field of weight $i$, i.e. $[\zeta, v_i] = iv_i$. The Lie derivative $L_{v_i}\zeta$ can be computed as 
$$ L_{v_i}\zeta = [v_i, \zeta] = -iv_i.$$ 
For the Poisson bivector $\theta$, we have $L_{(\psi'_{n+1})_*\zeta_1} \theta = L_{\zeta_2} \theta = -l\cdot \theta$; and hence $L_{v_i}\theta = 0$.  
  
%By the symplectic form $\omega$, it can be identified with the de Rham complex 
%$$ \Gamma(X_{reg}, O_{X_{reg}})(i + l) \stackrel{d}\to \Gamma (X_{reg}, \Omega^1_{X_{reg}})(i+l)  
%\stackrel{d}\to \Gamma (X_{reg}, \Omega^2_{X_{reg}})(i + l) $$ 
%By a similar argument to the proof of Proposition 3.2 of [Na 2] and Remark 3.4 of [ibid], 
%we see that these complexes are exact. Since $L_{\zeta_1}\omega = L_{\zeta_2}\omega 
%= l\cdot \omega$, we have $\delta_1(v_i) = 0$. By the exact sequences above, there is an 
%element $f_{i+l} \in \Gamma (X_{reg}, O_{X_{reg}})(i + l)$ such that $\delta_0 (f_{i+l}) = v_i$. 
%Note that $i(\delta_0(f_{i+l}))\omega = df_{i+l}$. We are going to prove that $L_{\delta_0(f_{i+l})}
%\zeta = - i\cdot \delta_0(f_{i+l})$. By the identification of $\Theta_X$ with $\Omega^1_X$ by 
%$\omega$ we identify $\zeta$ with a 1-form $\theta := i(\zeta)\omega$. As $i(\delta_0f_{i+l})\omega = df_{i+l}$, it is enough to show that $L_{\delta_0f_{i+l}} \theta = - i\cdot df_{i+l}$. This follows from Cartan's formula
%$$ L_{\delta_0f_{i+l}} \theta = d(i(\delta_0f_{i+l})\theta) + i(\delta_0f_{i+l})d\theta.$$ 
%In fact we have $$d\theta = d (i(\zeta)\omega) = L_{\zeta}\omega = l\cdot \omega,$$
%and $$i(\delta_0f_{i+l})\theta = i(\delta_0f_{i+l})i(\zeta)\omega = 
%- i(\zeta)i(\delta_0f_{i+l})\omega = -i(\zeta)df_{i+l} = -(i+l)f_{i+l}.$$
 
By this observation, if we put $$\psi_{n+1} := \psi'_{n+1} + t^{n+1}\Sigma_{i \ne 0}(1/i)v_i,$$ then $\psi_{n+1}$ is still a Poisson automorphism of $X_{n+1}$ (because $L_{v_i}\theta = 0$) and one can write 
$$(\psi_{n+1})_*\zeta_1 - \zeta_2 = t^{n+1}v_0.$$ 
Finally we show that $(\psi_{n+1})_*\zeta_1 = \zeta_2$ by using the fact that these are both 
integrated to $\mathbf{C}^*$-actions.        
Consider the two $\mathbf{C}^*$-actions on $X_{n+1}$ generated by  $(\psi_{n+1})_*\zeta_1$ and $\zeta_2$. As both vector fields are $\mathbf{C}^*$-invariant 
(with respect to any one of the two $\mathbf{C}^*$-actions), these $\mathbf{C}^*$ actions mutually commute. 
We now prove that these $\mathbf{C}^*$-actions are the same by the induction on $n$ (the 
index of $T_n$). The coordinate ring $\mathcal{R}$ of $X_{n+1}$ is 
isomorphic to $R \otimes_{\mathbf C}\mathbf{C}[t]/(t^{n+2})$, where $X = \mathrm{Spec}(R)$. 
We may assume that one of the $\mathbf{C}^*$-actions corresponds to the usual grading  
$\oplus_{i \geq 0} (R_i \oplus tR_i \oplus ... \oplus t^{n+1}R_i)$. Let us consider the weight $i$ 
eigenspace $V_i$ of another $\mathbf{C}^*$-action. Since two $\mathbf{C}^*$-actions 
are compatible, $V_i$ decomposes as $V_i = \oplus_j V_{i,j}$ where $V_{i,j}$ is a subspace 
of  $R_j \oplus tR_j \oplus ... \oplus t^{n+1}R_j$.   
By the induction hypothesis, we have $V_{i,j} \subset t^{n+1}R_j$ if $j \ne i$. 
But, the weight $i$ eigenspace of $t^{n+1}R $ with respect to the 2-nd $\mathbf{C}^*$-action 
also coincides with $t^{n+1}R_i$ 
because $t^{n+1}R  = t^{n+1}\mathcal{R} = (t^{n+1})\otimes_{\mathbf C} R$.
This means that $V_{i,j} = 0$ if $j \ne i$. Therefore $V_i = V_{i,i} \subset 
R_i \oplus tR_i \oplus ... \oplus t^{n+1}R_i$ and we conclude that  
the two $\mathbf{C}^*$-actions are the same.  
   
Now, the composite   
$$\{X_n\}_{n \geq 0} \stackrel{\beta_n \circ \psi_n}\to \{X \times T_n\}_{n \geq 0}$$ is a 
$\mathbf{C}^*$-equivariant Poisson isomorphism. 
By using the $\mathbf{C}^*$-actions 
of both sides, we then get a desired $\mathbf{C}^*$-equivariant Poisson isomorphism 
$\mathcal{X} \times_T \hat{T} \cong  X \times \hat{T}$ over $\mathrm{Spec}(\hat{A})$. 
Q.E.D.  

\begin{Prop} \label{(3.3)} Under the same assumption of \ref{(3.2)}, all fibres are 
isomorphic as conical symplectic varieties. 
\end{Prop}

{\em Proof}. Let us consider the two $T$-schemes $\mathcal{X}$ and $X \times T$ with 
$\mathbf{C}^*$-actions. We define a functor $$\mathrm{Hom}^{\mathbf{C}^*}_T(\mathcal{X}, X \times T): (T-schemes) \to (Set)$$ by 
$T' \to \mathrm{Hom}_{T'}^{\mathbf{C}^*}(\mathcal{X} \times_T T', X \times T')$. 
Then it is a functor of locally finite presentation. By Artin's approximation theorem [Ar], if we are given a $\mathbf{C}^*$-equivariant morphism $f: \mathcal{X} \times_T \hat{T} \to X \times \hat{T}$, then there is a pointed algebraic scheme $s_0 \in S$ (i.e. a pointed scheme of finite type over $\mathbf{C}$) together with an \'{e}tale map $h: (S, s_0) \to (T,0)$ and a $\mathbf{C}^*$-equivariant morphism $g: \mathcal{X} \times_T S \to X \times 
S$ such that $g(s_0): \mathcal{X} \times_T k(s_0) \to X \times k(s_0)$ coincides with $f(0): 
\mathcal{X} \times_T k(0) \to X \times k(0)$. 
We apply this to the morphism $f$ in 
\ref{(3.2)}.  As $f$ is an isomorphism, we may assume that $g$ is also an isomorphism, if necessary, by shrinking $S$ suitably.
Then we have a $\mathbf{C}^*$-equivariant isomorphism 
$\mathcal{X} \times_T S \cong X \times S$.  This implies that $\mathcal{X}_{h(s)}$ is 
isomorphic to $X$ as a $\mathbf{C}^*$-varieties for any closed point $s \in S$. By Theorem 3.1 of \cite{Na 2}, 
two conical symplectic varieties having the symplectic 2-forms of the same weight 
are isomorphic if they are isomorphic as $\mathbf{C}^*$-varieties.  Q.E.D.

Now, by \ref{(2.10)} and \ref{(3.3)}, one has

\begin{Thm} \label{(3.4)} For positive integers $N$ and $d$, there are only finite number of 
conical symplectic varieties of dimension $2d$ with maximal weights $N$, up to an isomorphism \end{Thm}

Department of Mathematics, Faculty of Science, Kyoto University 

e-mail address: namikawa@math.kyoto-u.ac.jp

\end{document}